\documentstyle[12pt]{amsart}
\newtheorem{lemma}{Lemma}[section]
\newtheorem{proposition}[lemma]{Proposition}

\newtheorem{theoremi}{Theorem}

\newtheorem{problem}{Problem}
\newtheorem{corollaryi}{Corollary}

\newcommand{\bC}{{\Bbb C}}
\newcommand{\bP}{{\Bbb P}}
\newcommand{\bQ}{{\Bbb Q}}
\newcommand{\bR}{{\Bbb R}}
\newcommand{\cA}{{\cal A}}
\newcommand{\cC}{{\cal C}}
\newcommand{\cH}{{\cal H}}

\newcommand{\cM}{{\cal M}}
\newcommand{\cN}{{\cal N}}
\newcommand{\cO}{{\cal O}}
\newcommand{\cU}{{\cal U}}
\newcommand{\cW}{{\cal W}}

\newcommand{\cS}{{\cal S}}
\newcommand{\cY}{{\cal Y}}

\newcommand{\orho}{\overline{\rho}}
\newcommand{\oS}{\overline{\cal S}}
\newcommand{\ooS}{\overline{\overline{\cal S}}}

\newcommand{\om}{\overline{m}}
\newcommand{\tbP}{\widetilde{\bP}}
\newcommand{\tM}{\widetilde{\cal M}}

\newcommand{\lra}{\longrightarrow}
\newcommand{\inj}{\hookrightarrow}
\newcommand{\surj}{\longrightarrow}

\begin{document}
\title[Density and completeness]{Density and completeness of
subvarieties of moduli spaces of curves or abelian varieties}
\author {E. Izadi}
\date{July 1996}
\address{Department of Mathematics, Boyd Graduate Studies Research Center, University of Georgia, Athens, GA 30602-7403}
\email{izadi@@math.uga.edu}
\thanks{The author was partially supported by an NSF grant.} 

\begin{abstract}
Let $V$ be a subvariety of codimension $\leq g$ of the moduli space $\cA_g$
of principally polarized abelian varieties of dimension $g$ or of the
moduli space $\tM_g$ of curves of compact type of genus $g$. We prove that
the set $E_1(V)$ of elements of $V$ which map onto an elliptic curve is
analytically dense in $V$. From this we deduce that if $V \subset \cA_g$ is
complete, then $V$ has codimension equal to $g$ and the set of elements of
$V$ isogenous to a product of $g$ elliptic curves is countable and
analytically dense in $V$. We also prove a technical property of the
conormal sheaf of $V$ if $V \subset \tM_g$ (or $\cA_g$) is complete of
codimension $g$.
\end{abstract}
\maketitle

\vskip25pt
\begin{center}
{\large \sc Introduction}
\end{center}
\vskip20pt

Let $\cM_g$ be the moduli space of smooth curves of genus $g \geq 2$
and let $\cA_g$ be the moduli space of principally polarized abelian
varieties ({\em ppav}) of dimension $g$ over $\bC$. A
(Deligne-Mumford) stable curve of genus $g$ is a reduced, connected
and complete curve of arithmetic genus $g$ with only nodes as
singularities and with finite automorphism group. We say that a stable
curve is of compact type if its generalized jacobian is an abelian
variety. We denote by $\tM_g$ the moduli space of stable curves of
compact type and genus $g$ over $\bC$. By ``density'' we always mean
``analytic density'' unless we specify otherwise.

Given a subvariety $V$ of $\cM_g$ or $\tM_g$ and an integer $q$ between $1$
and $g/2$, let $E_q(V)$ be the subset of $V$ parametrizing curves whose
jacobian contains an abelian variety of dimension $q$. We define $E_q(V)$
for $V$ a subvariety of $\cA_g$ in a similar fashion. It is well-known that
$E_q(\cA_g)$ is dense in $\cA_g$ for all $q$. Colombo and
Pirola pose the following question in \cite{colpi}
\begin{problem} 
When is $E_q(V)$ dense in $V$?
\end{problem}

Colombo and Pirola give a sufficient condition for the density of
$E_q(V)$ in $V$. They then show that $E_1(V)$ is dense in $V$ for all
subvarieties $V$ of $\cM_g$ of codimension at most $g-1$. They deduce from
this a second proof of the noncompleteness of codimension $g-1$
subvarieties of $\tM_g$ which was originally proved by Diaz in
\cite{diaz2}, Corollary page 80 (Colombo and Pirola prove the
noncompleteness of codimension $g-1$ subvarieties of $\tM_g$ which meet
$\cM_g$ ; however, if the subvariety is contained in $\tM_g \setminus
\cM_g$, its noncompleteness can be easily seen by mapping it (or, in some
cases, a double cover of it) to a moduli space of curves of lower genus).

Using the condition of Colombo and Pirola, we show
\begin{theoremi}
Suppose that $g \geq 2$. Let $V$ be a subvariety of codimension at most $g$
of $\cM_g$ or $\cA_g$, then $E_1(V)$ is dense in $V$.
\end{theoremi}
This result brings out another fundamental difference
between the moduli spaces in characteristic zero and in positive
characteristic (see Section \ref{sectpb} below).

We also obtain
\begin{corollaryi}
Suppose that $g \geq 2$. Let $V$ be a complete subvariety of codimension at
most $g$ of $\cA_g$. For all $i \in \{ 1, ..., g \}$, denote by
$E_{1,i}(V)$ the subset of $V$ parametrizing ppav's isogenous to a
product of a ppav of dimension $g-i$ and $i$ elliptic curves. Then
\begin{enumerate}
\item the variety $V$ has codimension exactly $g$,
\item for all $i \in \{ 1, ..., g \}$, any irreducible component $Z$ of
$E_{1,i}(V)$ has the expected dimension $\frac{(g-i)(g-i+1)}{2} + i
-g$; furthermore, the variety $Z$ parametrizes ppav's isogenous to a
product of $i$ {\em fixed} elliptic curves (depending only on $Z$) and
some ppav of dimension $g-i$,
\item for all $q, 1 \leq q \leq g/2$, any irreducible component of
$E_q(V)$ has the expected dimension $\frac{(g-i)(g-i+1)}{2} +
\frac{i(i+1)}{2} -g$,
\item the set $E_{1,g}(V)$ is dense in $V$. In particular, the set
$E_{1,g}(V)$ is (countable) infinite and, for all $i \in \{1, ..., g \}$,
for all $q \in \{1, ..., g/2 \}$, the sets $E_{1,i}(V)$ and $E_q(V)$ are
dense in $V$ (since they contain $E_{1,g}(V)$).
\end{enumerate}
\end{corollaryi}

An immediate consequence of Corollary $1$ is that
complete subvarieties of $\cA_g$ have codimension at least $g$. There
is also a proof of this last fact in \cite{oortcsubv} ($2.5.1$ page
$231$). There are no known examples of complete subvarieties of
codimension $g$ of $\cA_g$ (or $\tM_g$) except for $g=2$, although $g$
is the best known lower bound for the codimension of complete
subvarieties of $\cA_g$ (and $\tM_g$). It is conjectured in
\cite{oortcsubv} ($2.3$ page $230$) that for $g \geq 3$ the codimension of
a complete subvariety of $\cA_g$ is at least $g+1$.

Let $\cH_g$ be the locus of hyperelliptic curves in $\cM_g$. Let
$\tM_g'$ and $\cA_g'$ be respectively the moduli space of curves of
compact type with level $n$ structure (for some fixed $n \geq 3$) and
the moduli space of ppav's with level $n$ structure. Let $s_a : \cA_g'
\lra \cA_g$ and $s_c : \tM_g' \lra \tM_g$ be the natural morphisms. It
is well-known that $\cA_g'$ and $\cM_g' := s_c^{-1}(\cM_g)$ are smooth
and that there is a universal family of abelian varieties with level
$n$ structure on $\cA_g'$ and a universal family of (smooth) curves
with level $n$ structure on $\cM_g'$. By a ``universal'' family of
curves or abelian varieties with level $n$ structure we mean a family
which solves the moduli problem for curves or abelian varieties with
level $n$ structure. We note that the only properties of $\tM_g'$,
$\cM_g'$ and $\cA_g'$ we need are the smoothness of $\cM_g'$ and
$\cA_g'$ and the existence of the universal families (also note that
with non-abelian level structures or Prym-level structures, one
can get smooth covers of $\tM_g'$ as well (see \cite{looijenga1} and
\cite{djongpik})). We have the following technical consequence of our
results.
\begin{corollaryi}
Suppose that $g \geq 3$. Let $V$ be a complete subvariety of
codimension $g$ of $\tM_g$ or $\cA_g$. If $V \subset \cA_g$, let $V_0$
be the smooth locus of $s_a^{-1}(V)$. If $V \subset \tM_g$, let $V_0$
be the smooth locus of $s_c^{-1}(V \cap (\cM_g \setminus
\cH_g))$. Then the conormal bundle to $V_0$ is isomorphic to the
tensor product of the Hodge bundle (the pushforward of the sheaf of
relative one-forms on the universal abelian (or jacobian) variety)
with a subline bundle of the Hodge bundle.
\end{corollaryi}

Finally, we point out that, if $\cA_{g,d}$ denotes the moduli space of
abelian varieties of dimension $g$ and polarization type $d=(d_1,...,d_g)$,
then there is a finite correspondence between $\cA_g$ and $\cA_{g,d}$
so that our results remain valid if we replace $\cA_g$ with $\cA_{g,d}$.

\vskip10pt
{\sc Acknowledgments}
\vskip10pt

I wish to thank F. Oort for mentioning these questions to me
and also for stimulating discussions.

\vskip20pt

\begin{center}
{\sc Notation}
\end{center}

\vskip15pt

For any vector space or vector bundle (resp. affine cone) $E$, we
denote by $\bP(E)$ the projective space (resp. projective variety) of
lines in $E$ and by $E^*$ its dual vector space or vector bundle. We
let $E^{\otimes m}$, $S^mE$ and $\Lambda^mE$ respectively be the
$m$-th tensor power, the $m$-th symmetric power and the $m$-th
alternating power of $E$. For any linear map of vector spaces or
vector bundles $l : E \lra F$, we denote by $\overline{l} : {\Bbb
P}(E) \lra {\Bbb P}(F)$ its projectivization.

For any variety $X$ and any point $x \in X$, we denote by $T_xX$ the
Zariski tangent space to $X$ at $x$ and by $T_x^*X$ the dual of
$T_xX$. We denote by $X_{sm}$ the subvariety of smooth points of $X$.

For a ppav $A$, we let $\rho : H^0(\Omega^1_A)^{\otimes 2} \lra
S^2H^0(\Omega^1_A)$ be the natural linear map with kernel
$\Lambda^2H^0(\Omega^1_A)$. For a subvariety $V$ of $\cA_g'$ the
restriction to $V$ of the universal family on $\cA_g'$ gives a family
$\cA_V$ of ppav's on $V$ (we forget the level $n$ structure). For a point
$t$ of $V$, we let $A_t$ be the fiber of $\cA_V$ at $t$. The
Zariski-tangent space to $\cA_g'$ at $t$ can be canonically identified with
$S^2H^0(\Omega^1_{A_t})^*$. We denote by $\pi_a : S^2H^0(\Omega^1_{A_t})
\lra T^*_tV$ the codifferential at $t$ of the embedding $V \inj \cA_g'$.

For a smooth curve $C$, we denote by $\omega_C$ the canonical sheaf of $C$
and let $\kappa C$ be the image of $C$ in the dual projective space
$|\omega_C|^*$ of the linear system $|\omega_C|$ by the natural morphism
associated to this linear system. If $A= JC$ is the jacobian of a smooth
curve $C$, then $H^0(\Omega^1_A) \cong H^0(\omega_C)$. Let $m :
S^2H^0(\omega_C) \lra H^0( \omega_C^{\otimes 2})$ be multiplication
and put $\mu := m \rho$. For $V
\subset \cM_g'$ the restriction to $V$ of the universal family on $\cM_g'$
gives a family $\cC_V$ of curves on $V$ (again we forget the level $n$
structure). For $t \in V$ we let $C_t$ be the fiber of $\cC_V$ at $t$. The
Zariski-tangent space to $\cM_g'$ at $t$ can be canonically identified with
$H^0(\omega_{C_t}^{\otimes 2})^*$. We let $\pi : H^0(\omega_{C_t}^{\otimes
2}) \lra T^*_tV$ be the codifferential at $t$ of the embedding $V \inj
\cM_g'$.

For a stable curve $C$ of compact type (resp. a ppav $A$), we will call the
corresponding point of $\tM_g$ (resp. $\cA_g$) the moduli point of $C$
(resp. $A$).

\section{The proofs}
\label{sectpf}
In this section we give the proof of Theorem 1 and its corollaries.

We first consider the case where $V$ is contained in $\cM_g$. We may
and will replace $V$ with its inverse image in $\cM_g'$.

The relative jacobian of $\cC_V$ gives us a family of ppav's on $V$.
We can therefore apply Theorem (1) on page 162 of \cite{colpi}: to
show that $E_1(V)$ is dense in $V$ it is enough to prove the
following:

There exists a Zariski dense (Zariski-)open subset $U$ of $V$,
contained in the smooth locus $V_{sm}$ of $V$, such that, for all $t \in
U$, there is a subvector space $W$ of $H^0(\omega_{C_t})$ which has
dimension $1$ and is such that the composition
\[W \otimes W^{\perp} \stackrel{\mu}{\inj} H^0(\omega_{C_t}^{\otimes 2}) \stackrel{\pi}{\surj}
T^*_tV
\]
is injective. Here $W^{\perp}$ is the orthogonal complement of $W$
with respect to the hermitian form on $H^0(\omega_{C_t})$ induced by the
natural polarization of $JC_t$. We sketch briefly how this condition is
obtained in the more general case where $W$ has dimension $q$ with $1
\leq q \leq g/2$ and $V \subset \cA_g'$ has any dimension $\geq
q(g-q)$.

To prove the density of $E_q(V)$ in $V$, it is enough to show that there is
a Zariski dense open subset $U$ of $V_{sm}$, such that,
for all $t \in U$, there is an analytic neighborhood $U'$ of $t$, $U'
\subset U$, such that $E_q(V) \cap U'$ is dense in $U'$. An abelian
variety $A$ contains an abelian subvariety of dimension $q$ if and
only if $H^0(\Omega^1_{A})$ contains a $q$-dimensional ${\Bbb
C}$-subvector space which is the tensor product with ${\Bbb R}$ of a
vector subspace of dimension $2q$ of $H^1(A, {\Bbb Q})$ (after
identifying $H^0(\Omega^1_{A})$ with $H^1(A, \bR) \cong H^1(A, \bQ)
\otimes \bR$ as real vector spaces). Let $t$ be an element of
$V_{sm}$. For a contractible analytically open set $U' \ni t$
contained in $V_{sm}$, let $F_{U'}$ be the Hodge bundle over
$U'$. Then one can trivialize $F_{U'}$ as a real vector
bundle. Therefore the grassmannian bundle of $2q$-dimensional real
subvector spaces of $F_{U'}$ is isomorphic to $U' \times G_{\Bbb
R}(2q,2g)$, where $G_{\Bbb R}(2q,2g)$ is the Grassmannian of
$2q$-dimensional ${\Bbb R}$-subvector spaces of $H^1(A_t,
\bR)$. Hence there is a well-defined map $\Phi : G(q,F_{U'}) \lra G_{\Bbb
R}(2q,2g)$ where $G(q,F_{U'})$ is the Grassmannian of $q$-dimensional
${\Bbb C}$-subvector spaces of $F_{U'}$: The map $\Phi$ sends a
$q$-dimensional complex subvector space of $H^0(\Omega^1_{A_s})$ (with $s
\in U'$) to the image of its underlying real vector space under the
isomorphism $H^1(A_s, \bR) \stackrel{\cong}{\lra} H^1(A_t, \bR)$
obtained from the $\bR$-trivialization of $F_{U'}$. Let
$G_{\bQ}(2q,2g) \subset G_{\bR}(2q,2g)$ be the Grassmannian of
$2q$-dimensional $\bQ$-subvector spaces of $H^1(A_t, \bQ)$ and let $p
: G(q,F_{U'}) \lra U'$ be the natural morphism. Then $E_q(V) \cap U' =
p(\Phi^{-1}(G_{\bQ}(2q,2g)))$. To prove the density of $E_q(V) \cap
U'$ in $U'$, it is enough to prove that there is a subset $\cY$ of
$G(q,F_{U'})$ such that $p(\cY) = U'$ and $\Phi^{-1}(G_{\bQ}(2q,2g))
\cap \cY$ is dense in $\cY$. Since $G_{\bQ}(2q,2g)$ is dense in
$G_{\bR}(2q,2g)$, it is enough to find $\cY$ such that $p(\cY) = U'$
and $\Phi |_{\cY}$ is an open map. If $\Phi$ has maximal rank (i.e.,
the differential $d \Phi$ of $\Phi$ is surjective) everywhere on
$\cY$, then $\Phi |_{\cY}$ is an open map. Therefore $E_q(V) \cap U'$
is dense in $U'$ if for every $s \in U'$ there is a $q$-dimensional
$\bC$-subvector space $W$ of $H^0(\Omega^1_{A_s})$ such that $d \Phi$
is surjective at $(W,s) \in G(q,F_{U'})$ (then $\cY$ would be the set
of such $(W,s)$). The tangent space $T_{(W,s)}G(q,F_{U'})$ is
isomorphic to $W
\otimes W^{\perp} \oplus T_sU'$, the tangent space to $G_{\Bbb R}(2q,2g)$
at $\Phi(W,s)$ is isomorphic to $W \otimes W^{\perp} \oplus \overline{W
\otimes W^{\perp}} \cong W \otimes W^{\perp} \oplus (W \otimes
W^{\perp})^*$ and the restriction of $d \Phi$ to the $W \otimes
W^{\perp}$ summand of $T_{(W,s)}G(q,F_{U'})$ is an isomorphism onto the $W
\otimes W^{\perp}$ summand of
$T_{\Phi(W,s)}G_{\Bbb R}(2q,2g)$. Therefore $d \Phi$ is surjective if
and only if the map it induces $T_sU' = \frac{T_{(W,s)}G(q,F_{U'})}{W
\otimes W^{\perp}} \lra (W \otimes W^{\perp})^* =
\frac{T_{\Phi(W,s)}G_{\Bbb R}(2q,2g)}{W \otimes W^{\perp}}$ is
surjective, i.e., if and only if the dualized map $W \otimes W^{\perp}
\lra T^*_sU'$ is injective. Let $F$ be the Hodge bundle over the
Siegel upper half space $\cU_g$. The inclusion $U' \inj \cA_g'$ lifts
to an inclusion $U' \inj \cU_g$ because $U'$ is contractible and there
is a family of ppav's on $U'$ (the restriction of $\cA_V$). Factoring
$\Phi$ through the Grassmannian of $q$-planes in $F$ over $\cU_g$, the
map $W \otimes W^{\perp} \lra T_s^*U'$ can be seen to be the composition
\[ W \otimes W^{\perp} \stackrel{\rho}{\lra} S^2H^0(\Omega^1_{A_s})
\stackrel{\pi_a}{\lra} T_s^*U' = T_s^*V \: .
\]
For $V$ contained in $\cM_g'$, we have $\pi_a = \pi m$.

\vskip10pt

Clearly, if $\pi \mu : W \otimes H^0(\omega_{C_t}) \lra T^*_tV$ is
injective, then so is $\pi \mu : W \otimes W^{\perp} \lra T^*_tV$. In view
of this (and also for use in the proof of Corollary 2) we show:

\begin{proposition}
Suppose that $g \geq 3$. Let $V$ be a subvariety of codimension at most $g$
of $\cM_g'$. Let $t$ be a point of $V_{sm}$ and let $N$ be the kernel of
$\pi : H^0(\omega_{C_t}^{\otimes 2}) \surj T^*_tV$.
\begin{enumerate} 
\item Suppose that $C_t$ is non-hyperelliptic. Suppose that, for any
one-dimensional subvector space $W$ of $H^0(\omega_{C_t})$, the map $\pi
\mu : W \otimes H^0(\omega_{C_t}) \lra T^*_tV$ is {\em not} injective. Then
$V$ has codimension exactly $g$ and there is a one-dimensional subvector
space $W_N$ of $H^0(\omega_{C_t})$ such that $N = \mu (W_N \otimes
H^0(\omega_{C_t}))$.
\item Suppose that $C_t$ is hyperelliptic and that $V$ is {\em not}
transverse to $\cH_g' := s_c^{-1}(\cH_g)$ at $t$ (i.e., the sum $T_tV + T_t
\cH_g' \subset T_t \cM_g'$ is {\em not} equal to $T_t \cM_g'$). Then there
exists a one-dimensional subvector space $W$ of $H^0(\omega_{C_t})$ such
that the map $\pi \mu : W \otimes H^0(\omega_{C_t}) \lra T^*_tV$ is
injective.
\end{enumerate}
\label{propNcurve}
\end{proposition}
{\em Proof :} Consider the composition
\[ {\Bbb P}(H^0(\omega_{C_t})^{\otimes 2}) \stackrel{\orho}{\lra}
{\Bbb P}(S^2H^0(\omega_{C_t})) \stackrel{\overline{m}}{\lra} {\Bbb
P}(H^0(\omega_{C_t}^{\otimes 2})) \: .
\]
The kernel of $m$ is the space $I_2(C_t)$ of quadratic forms vanishing on
$\kappa C_t$. Hence the {\em rational} map $\om$ is the projection with
center $\bP (I_2(C_t))$. Let $\oS$ be the image by $\orho$ of the Segre
embedding $\cS$ of ${\Bbb P}(H^0(\omega_{C_t})) \times {\Bbb
P}(H^0(\omega_{C_t}))$ in ${\Bbb P}(H^0(\omega_{C_t})^{\otimes 2})$. Let
$N'$ be the set of rank $2$ symmetric tensors in $S^2H^0(\omega_{C_t})$
which lie in $m^{-1}(N)$ (then $\bP(N')$ is the {\em reduced} intersection
of $\oS$ and $\overline{m}^{-1}(\bP(N))$).

Suppose that for all $W \subset H^0(\omega_{C_t})$ of dimension $1$,
the map $\mu : W \otimes H^0(\omega_{C_t}) \lra T_t^*V$ is not
injective, i.e., for all $w \in H^0(\omega_{C_t})$, there is $w' \in
H^0(\omega_{C_t})$ such that $\mu (w \otimes w') \in N$. This implies that
the dimension of $\bP(N')$ is at least $g-1$. We will show below that this
does not happen if $C_t$ is hyperelliptic and $V$ is not transverse to
$\cH_g'$ at $t$. If $C_t$ is non-hyperelliptic, we will show that this
implies that $\bP(N')$ is a linear subspace of $\oS$ and that its inverse
image in $\cS$ is the union of two linear subspaces of $\cS$ which are two
fibers of the two projections of $\cS$ onto $\bP^{g-1}$ and are exchanged
under the involution of $\cS$ which interchanges the two
$\bP^{g-1}$-factors of $\cS$. The proposition will then easily follow from
this.

Suppose first that $C_t$ is {\em non}-hyperelliptic. Then $m :
S^2H^0(\omega_{C_t}) \lra H^0(\omega_{C_t}^{\otimes 2})$ is onto (see
\cite{ACGH} page 117). We have
\begin{lemma}                                                     
Suppose $g= 2$ or $g \geq 3$ and $C_t$ is non-hyperelliptic. Suppose
that for all $W \subset H^0(\omega_{C_t})$ of dimension $1$, the map
$\mu : W \otimes H^0(\omega_{C_t}) \lra T_t^*V$ is not injective. Then the
map $\bP(N') \lra {\Bbb P}(N)$ is generically one-to-one.
\label{N'isPN}
\end{lemma}
{\em Proof :} If not, then, for all $w \in H^0(\omega_{C_t})$, there
exists $w', w_1, w_1'\in H^0(\omega_{C_t})$ such that $w w' := \rho(w
\otimes w')$ and $w_1 w_1' := \rho(w_1 \otimes w_1')$ are not
proportional but $m (w w')$ and $m (w_1 w_1')$ are proportional
elements of $N$. Therefore, supposing $w$ general, there exits
$\lambda \in \bC, \lambda \neq 0$, such that the element $\lambda
ww'-w_1w_1'$ of $S^2H^0(\omega_{C_t})$ lies in $I_2(C_t)$, i.e.,
defines a quadric $q(w)$ of rank $3$ or $4$ (in the canonical space $|
\omega_{C_t} |^*$) which contains $\kappa C_t$ (the
canonical curve $\kappa C_t$ is not contained in any quadric of rank
$\leq 2$ since it is nondegenerate). If $g \leq 3$, this is impossible
because in that case $I_2(C_t)=0$. If $g \geq 4$, the intersection $L$
of the two hyperplanes in $|\omega_{C_t}|^*$ with equations $w$ and
$w_1$ is an element of a ruling of the quadric $q(w)$. Therefore $L$
cuts a divisor of a $g^1_d$ (a $g^1_d$ is a pencil of divisors of
degree $d$) on $C_t$ with $d \leq g-1$ (see \cite{AM}, Lemmas 2 and 3
page 192). Therefore the divisor of zeros of $w$ on $C_t$ contains a
divisor of a $g^1_d$. By the uniform position Theorem (see \cite{ACGH}
Chapter $3, \S 1$) this does not happen for $w$ in some nonempty
Zariski-open subset of $H^0(\omega_{C_t}) \setminus \{ 0 \}$.  \hfill
\qed

\vskip20pt
Therefore, since the dimension of $\bP(N')$ is at least $g-1$ and the
dimension of ${\Bbb P}(N)$ is at most $g-1$, the map $\bP(N') \lra \bP (N)$
is birational and $\bP(N')$ and ${\Bbb P}(N)$ have
both dimension $g-1$.

{\em This proves, in particular, that $V$ has codimension {\em
exactly} $g$.}

Since no quadrics of rank $\leq 2$ contain $\kappa C_t$, the center
$\bP I_2(C_t)$ of the projection $\om$ does not intersect $\oS$. In
particular, the space $\bP I_2(C_t)$ does not intersect
$\bP(N')$. Therefore $\om$ restricts to a birational {\em morphism}
$\bP(N') \surj \bP(N)$ and, since ${\Bbb P}(N)$ is a linear subspace
of ${\Bbb P}(H^0(\omega_{C_t}^{\otimes 2}))$, the degree of $\bP(N')$
(in the projective space $\bP (S^2H^0(\omega_{C_t}))$) is equal to the
(generic) degree of the map $\bP(N') \surj \bP(N)$. Hence $\bP(N')$ is
a linear subspace of ${\Bbb P}(S^2H^0(\omega_{C_t}))$ and
$\overline{m}$ restricts to an {\em isomorphism} $\bP(N')
\stackrel{\cong}{\lra} {\Bbb P}(N)$.

Let $N''$ be the cone of decomposable tensors in
$H^0(\omega_{C_t})^{\otimes 2}$ which lie in $\mu^{-1}(N)$ (then
$\bP(N'')$ is the {\em reduced} inverse image of $\bP(N')$ in $\cS
\subset {\Bbb P}(H^0(\omega_{C_t})^{\otimes 2})$). The map $\cS \lra
\oS$ is a finite {\em morphism} of degree $2$ ramified on the
diagonal. Therefore the map $\bP(N'') \lra \bP(N')$ is a morphism of
degree $\leq 2$. Since the diagonal of $\oS \cong S^2 {\Bbb
P}(H^0(\omega_{C_t}))$ is irreducible of dimension $g-1$ and spans
${\Bbb P}(S^2H^0(\omega_{C_t}))$, the space $\bP(N')$ intersects
this diagonal in a subvariety of dimension at most $g-2$. Therefore
the morphism $\bP(N'') \lra \bP(N')$ has degree $2$ and $\bP(N'')$ has
degree $2$ in $\bP(H^0(\omega_{C_t})^{\otimes 2})$.

If $\bP(N'')$ is irreducible, it spans a linear subspace $\tbP$ of $\bP
(H^0(\omega_{C_t})^{\otimes 2})$ of dimension $g$. This implies that $\bP
(\Lambda^2H^0(\omega_{C_t}))$ intersects $\tbP$ in exactly one point. For
$w_1, w_2 \in H^0(\omega_{C_t})$, let $w_1', w_2'$ be such that $\mu (w_1
\otimes w_1'), \mu (w_2 \otimes w_2') \in N$. For $w_i$ general, $w_i'$ is
not proportional to $w_i$ since $\bP(N'')$ intersects the diagonal of $\cS$
in a subvariety of dimension at most $g-2$. Therefore the lines spanned by
$w_1 \otimes w_1' - w_1' \otimes w_1$ and $w_2 \otimes w_2' - w_2' \otimes
w_2$ give us elements of $\tbP \cap {\Bbb P}(\Lambda^2H^0(\omega_{C_t}))$
which is a point. Therefore, for all $w_1, w_2 \in H^0(\omega_{C_t})$
general there exists $\lambda \in {\Bbb C}, \lambda \neq 0$, such that
\[ w_1 \otimes w_1' - w_1' \otimes w_1 = \lambda (w_2
\otimes w_2' -w_2' \otimes w_2)
\]
Complete $\{ w_1, w_2 \}$ to a general basis $\{ w_1, w_2, w_3, ..., w_g
\}$ of $H^0(\omega_{C_t})$ and write $w_i' = \sum_{1 \leq j \leq g}
a_{ij}w_j$ for $i = 1$ or $2$. Then from the equation above we deduce
$a_{1j}=a_{2j}=0$ for $j > 2$. Therefore, $w_1'$ belongs to the span of
$w_1$ and $w_2$. Repeating this argument with $w_1$ and $w_3$ instead
of $w_1$ and $w_2$, we see that $w_1'$ also belongs to the span of
$w_1$ and $w_3$. Hence $w_1'$ is proportional to $w_1$ (this is the
only part in the proof of Proposition \ref{propNcurve} where we need $g
\geq 3$). Contradiction.

Therefore $\bP(N'')$ is reducible, i.e., it is the union of two linear
subspaces of dimension $g-1$. We have
\begin{lemma}
Suppose $g \geq 2$. All linear subspaces of dimension $g-1$ of $\cS
\subset \bP ( H^0(\omega_{C_t})^{\otimes 2}) \cong {\Bbb P}^{g^2 -1}$
are elements of one of the two rulings of $\cS \cong \bP
(H^0(\omega_{C_t})) \times \bP (H^0(\omega_{C_t})) \cong {\Bbb
P}^{g-1} \times {\Bbb P}^{g-1}$.
\end{lemma}
{\em Proof :} Let $T$ be a linear subspace of dimension $g-1$ of $\cS$. Let
$p_1$ and $p_2$ be the two projections of $\cS \cong {\Bbb P}^{g-1} \times
{\Bbb P}^{g-1}$ onto its two factors. Let $H_i$ be a general element of
$p_i^* |{\cal O}_{{\Bbb P}^{g-1}}(1)|$ for $i= 1$ or $2$. Then $H_1 \cap T
\neq H_2 \cap T$ and $H_i$ does not contain $T$. In particular, the
intersection $H_i \cap T$ is either empty or of dimension $g-2$. The
divisor $H_1 \cup H_2$ is the intersection of a hyperplane $H$ in
${\Bbb P}^{g^2 -1}$ with $\cS$. Since $T$ is not contained in $H_1$
nor $H_2$, the hyperplane $H$ does not contain $T$ and hence $T \cap
H$ is a linear space of dimension $g-2$. Since the two intersections
$T \cap H_1 \neq T \cap H_2$ are both contained in the
$(g-2)$-dimensional linear space $T \cap H$ and are either empty or
have dimension $g-2$, we have either $H_1 \cap T =
\emptyset$ or $H_2 \cap T = \emptyset$. Suppose, for instance, that
$H_1 \cap T = \emptyset$. It is easily seen that $p_1^{-1}(p_1(H_1)) =
H_1$ implies $p_1(H_1) \cap p_1(T) = p_1(H_1 \cap T)$. Therefore
$p_1(T)$ does not intersect $p_1(H_1)$ which is a hyperplane in
$\bP^{g-1}$. Hence $p_1(T)$ is a point and $T$ is a fiber of $p_1$.
\hfill \qed

\bigskip

We deduce from the above Lemma that $\bP(N'') = \bP(N_1) \cup \bP(N_2)$
where $\bP(N_1)$ and $\bP(N_2)$ are elements of the two rulings of $\cS
\cong {\Bbb P}^{g-1} \times {\Bbb P}^{g-1}$. The spaces $\bP(N_1)$ and
$\bP(N_2)$ are exchanged by the involution which exchanges the two factors
of $\cS$ because $\bP(N'')$ is the inverse image of a linear subspace in
$\oS \cong S^2 {\Bbb P}^{g-1}$. Therefore there exists a one-dimensional
subvector space $W_N$ of $H^0(\omega_{C_t})$ such that, for instance, $N_1
= W_N \otimes H^0(\omega_{C_t})$ and $N_2 = H^0(\omega_{C_t}) \otimes
W_N$. So $N = \mu(N_1) = \mu (W_N \otimes H^0(\omega_{C_t})$. This proves
the Proposition in the non-hyperelliptic case.\\

Now suppose that $C_t$ is hyperelliptic and that $V$ is not transverse to
$\cH_g'$ at $t$, i.e., the subspaces $T_tV$ and $T_t \cH_g'$ do {\em not}
span $T_t \cM_g'$. Let $\iota$ be the hyperelliptic involution of
$C_t$. Let $H^0(\omega_{C_t}^{\otimes 2})^+$ and $H^0(\omega_{C_t}^{\otimes
2})^-$ be the subvector spaces of $H^0(\omega_{C_t}^{\otimes 2})$ of
$\iota$-invariant and $\iota$-anti-invariant quadratic differentials
respectively. Then $H^0(\omega_{C_t}^{\otimes 2})^+$ is the image of
$S^2H^0(\omega_{C_t})$ by $m$ and the conormal space to $\cH_g'$ at $t$ can
be canonically identified with $H^0(\omega_{C_t}^{\otimes 2})^-$. The
non-transversality of $V$ and $\cH_g'$ means that $N \cap
H^0(\omega_{C_t}^{\otimes 2})^- \neq \{ 0 \}$. This implies that $N$ is not
contained in $H^0(\omega_{C_t}^{\otimes 2})^+$. Since $N$ has dimension at
most $g$, the dimension of $N \cap H^0(\omega_{C_t}^{\otimes 2})^-$ is at
most $g-1$. Hence the dimension of $\bP(N \cap H^0(\omega_{C_t}^{\otimes
2})^+) = \bP(N) \cap \bP(H^0(\omega_{C_t}^{\otimes 2})^+) = \bP(N) \cap \om
(\bP(S^2H^0(\omega_{C_t})))$ is at most $g-2$. We have
\begin{lemma}
Suppose $g \geq 2$ and $C_t$ hyperelliptic. The map $\om : \oS
\surj \ooS := \om ( \oS )$ is a finite morphism of degree $ \frac{1}{2} \left(
\begin{array}{c} 2g-2 \\ g-1 \end{array} \right) $.
\label{oSooSfinite}
\end{lemma}
Note that the lemma finishes the proof of Proposition
\ref{propNcurve}: we saw above that the dimension of ${\Bbb P}(N) \cap \om
(\bP(S^2H^0(\omega_{C_t})))$ is at most $g-2$. A fortiori, since $\om
(\bP(S^2H^0(\omega_{C_t}))) \supset \ooS$, the dimension of $\bP(N)
\cap \ooS$ is at most $g-2$ and the dimension of $\bP(N')$ is at most
$g-2$ which is what we needed to show (see the paragraphs preceding
Lemma \ref{N'isPN}).\\

{\em Proof of lemma \ref{oSooSfinite}:} The map $\om : \oS \lra \ooS$
is a morphism if and only if the center $\bP(I_2(C_t))$ of the
projection $\om$ does not intersect $\oS$. This is the case because
the canonical curve $\kappa C_t$ is nondegenerate and hence not
contained in any quadrics of rank $\leq 2$.

Fix a nonzero element $w w' = \rho (w \otimes w')$ of
$S^2H^0(\omega_{C_t})$ and suppose that $w_1 w_1' \in S^2
H^0(\omega_{C_t})$ is not proportional to $w w'$ and $m (w_1 w_1') =
\lambda. m (w w')$ for some $\lambda \in \bC, \lambda \neq 0$. This is
equivalent to $Z(w) + Z(w') = Z(w_1) + Z(w_1')$ where $Z(w)$, for instance,
is the divisor of zeros of $w$ on the rational normal curve $\kappa C_t$. So
there are only a finite number of possibilities for $Z(w_1)$ and
$Z(w_1')$. This proves that $\om : \oS \lra \ooS$ is quasi-finite and hence
finite since it is proper. Any divisor of degree $g-1$ on $\kappa C_t \cong
\bP^1$ is the divisor of zeros of some element of $H^0(\omega_{C_t}) =
H^0(\cO_{\bP^1}(g-1))$, hence, since there are $ \frac{1}{2} \left(
\begin{array}{c} 2g-2 \\ g-1 \end{array} \right) $ ways to write a
fixed reduced divisor of degree $2g-2$ as a sum of two divisors of
degree $g-1$, the degree of $\om : \oS \lra \ooS$ is $ \frac{1}{2}
\left( \begin{array}{c} 2g-2 \\ g-1 \end{array}
\right) $. \hfill \qed

\vskip15pt
{\bf Proof of Theorem 1 in the case of curves:} As explained in the
beginning of this section, we need to find a Zariski-dense open subset $U$
of $V_{sm}$, such that, for all $t \in U$, there exists
$W \subset H^0(\omega_{C_t})$ ($W$ of dimension $1$) such that $\mu (W
\otimes W^{\perp}) \cap N = \{ 0 \}$.

First suppose $g \geq 3$. We may assume that $V$ is irreducible. If
$V$ is contained in $\cH_g'$, then $V$ is not transverse anywhere to
$\cH_g'$ and hence, by Proposition \ref{propNcurve}, we may take $U$
to be all of $V_{sm}$. If $V \not \subset \cH_g'$, take $U = V_{sm}
\setminus \cH_g'$. Suppose that there exists $t \in U$ such that, for
all $W \subset H^0(\omega_{C_t})$ of dimension $1$, we have $\mu (W
\otimes W^{\perp}) \cap N \neq \{ 0 \}$. Then, a fortiori, the
hypotheses of part 1 of Proposition \ref{propNcurve} are met and $N=
\mu (W_N \otimes H^0(\omega_{C_t}))$. Then every element of
$H^0(\omega_{C_t})$ is orthogonal to $W_N$. This is impossible given
that the hermitian form on $H^0(\omega_{C_t})$ is positive definite.

Now suppose $g=2$. Then $N$ has dimension $\leq 2$ and $\bP(N)$ has
dimension $\leq 1$. For each $W \subset H^0(\omega_{C_t})$ of dimension
$1$, the space $W^{\perp}$ also has dimension $1$ and hence $W \otimes
W^{\perp}$ has dimension $1$. The lines $W \otimes W^{\perp}$ form a
real analytic subset of $\bP (H^0(\omega_{C_t})^{\otimes 2})$ of real
dimension $2$. Since $\orho : \cS \lra \oS$ is finite, we deduce that
the lines $\rho(W \otimes W^{\perp}) = \mu (W \otimes W^{\perp})$
form a real analytic subset of $\bP (S^2H^0(\omega_{C_t})) =
\bP(H^0(\omega_{C_t}^{\otimes 2})) \cong \bP^2$ of real dimension $2$. An
easy computation (with coordinates) will show that this subset is not
contained in any projective line in $\bP(H^0(\omega_{C_t}^{\otimes 2}))$
and hence is not contained in $\bP(N)$. Hence there exists $W$ such that the
line $\mu (W \otimes W^{\perp})$ is not contained in $N$, in other
words $\mu(W \otimes W^{\perp}) \cap N = \{ 0 \}$. \hfill \qed

\vskip15pt
We now consider the case $V \subset \cA_g'$. As before, we first prove
\begin{proposition}
Suppose that $g \geq 3$. Let $V$ be a subvariety of codimension at most $g$
of $\cA_g'$. Let $t$ be a point of $V_{sm}$ and let $N$ be the kernel of
$\pi_a : S^2H^0(\Omega_{A_t}^1) \surj T^*_tV$. Suppose that, for any
one-dimensional subvector space $W$ of $H^0(\Omega_{A_t}^1)$, the map
$\pi_a \rho : W \otimes H^0(\Omega_{A_t}^1) \lra T^*_tV$ is {\em not}
injective. Then $V$ has codimension exactly $g$ and there is a
one-dimensional subvector space $W_N$ of $H^0(\Omega_{A_t}^1)$ such that $N
= \rho (W_N \otimes H^0(\Omega_{A_t}^1))$.
\label{propNppav}
\end{proposition}
{\em Proof :} If the map $\pi_a \rho : W \otimes H^0(\Omega_{A_t}^1) \lra
T^*_tV$ is {\em not} injective, then $\rho (W \otimes H^0(\Omega_{A_t}^1))
\cap N \neq \{ 0 \}$. If this holds for every $W \subset
H^0(\Omega_{A_t}^1)$ of dimension $1$, then ${\Bbb P}(N)$ has dimension
$g-1$ and is contained in $\oS \cong S^2{\Bbb P}(H^0(\Omega^1_A))$. It
follows that $V$ has codimension $g$. The rest of the argument is now
analogous to the proof of part $1$ of Proposition \ref{propNcurve} with $N'
= N$. \hfill \qed \vskip15pt

{\bf Proof of Theorem 1 in the case of abelian varieties:}
This proof is now as in the case of curves. \hfill \qed

\vskip15pt

{\bf Proof of Corollary 1:}

Let $V$ be a complete subvariety of codimension $g-d$ ($d \geq 0$) of
$\cA_g$. By Theorem 1, the set $E_1(V)$ is dense in $V$. In
particular, it is nonempty. Let $Y$ be an irreducible component of
$E_1(V)$. Let $r$ and $s$ be integers such that
for every ppav $A$ with moduli point in $Y$ there is an
elliptic curve $E$, a ppav $B$ and an isogeny $\nu : E \times B \lra
A$ of degree at most $r$ such that the inverse image of the principal
polarization of $A$ by $\nu$ is a polarization of degree at most
$s$. Let $Y'$ be an irreducible component of the variety parametrizing
such quadruples $(E,B,A, \nu)$. Then $Y'$ is a finite cover of $Y$.
The morphism $Y' \lra \cA_1$ which to $(E,B,A, \nu)$ associates the
isomorphism class of $E$ is constant since $Y'$ is complete (and
irreducible) and $\cA_1$ is affine.

For any irreducible component $Z$ of $E_1(\cA_g)$, there is a finite
correspondance between $Z$ and $\cA_{g-1} \times \cA_1$. In particular, the
codimension of $Z$ in $\cA_g$ is $\frac{g(g+1)}{2} - (\frac{g(g-1)}{2} +1)=
g-1$. The variety $Y$ is an irreducible component of the intersection of
$V$ with such a $Z$, hence there is a nonnegative integer $e_0$ such that
the codimension of $Y$ in $V$ is $g-1-e_0$. So the codimension of $Y$ in
$\cA_g$ is $g-d +g-1-e_0 = 2g-d-1-e_0$. Since $Y'$ maps to a point in
$\cA_1$, its image $V_1$ in $\cA_{g-1}$ by the second projection has
dimension equal to the dimension of $Y$. Therefore $V_1$ has dimension
$g(g+1)/2-(2g-d-1-e_0)=(g-1)g/2-(g-1-d-e_0)$, i.e., codimension $g-1-d-e_0
\leq g-1$ in $\cA_{g-1}$. By Theorem 1, the set $E_1(V_1)$ is dense in
$V_1$. In particular, the set $E_1(V_1)$ is nonempty. Let $Y_1$ be an
irreducible component of $E_1(V_1)$ and let $Y_1'$ be the analogue of
$Y'$ for $Y_1$. Then, as before, the variety $Y_1$ has codimension
$g-1-d-e_0+g-2-e_1$ in $\cA_{g-1}$ (for some nonnegative integer
$e_1$), the variety $Y_1'$ maps to a point in $\cA_1$ and its image
$V_2$ in $\cA_{g-2}$ has codimension $g-2-d-e_0-e_1$. Repeating the
argument, we obtain $V_i$ in $\cA_{g-i}$ of codimension $g-i-d-e_0 -
... - e_{i-1}$ containing $Y_i$ of codimension $g-i-d-e_0 -
... -e_{i-1}+g-i-1- e_i$ in $\cA_{g-i}$. For $i=g-2$, we can repeat
the argument one last time for $V_{g-2} \subset \cA_2$ to obtain
$Y'_{g-2}$ with image $V_{g-1}$ in $\cA_1$ with codimension
$1-d-e_0-...-e_{g-2}$. Since $\cA_1$ is affine, the variety $V_{g-1}$
is a point and $d=e_0=...=e_{g-2} = 0$. Therefore $Y$ has codimension
$2g-1$ in $\cA_g$, all the varieties $Y_i$ have codimension
$g-i+g-i-1=2g-2i-1$ in $\cA_{g-i}$, $V$ has codimension $g$ in $\cA_g$
and $V_i$ has codimension $g-i$ in $\cA_{g-i}$. In particular, the
first part of Corollary $1$ is proved.

For each $i$, there is an irreducible subvariety $Z_i$ of $V$ which
parametrizes ppav's isogenous to the product of an element of $V_i$
and $i$ fixed elliptic curves ($Z_1 = Y$) because all the maps $Y_i'
\lra \cA_1$ (and also $Y' \lra \cA_1$) are constant. It follows from
the above that $Z_i$ has the expected dimension
$\frac{(g-i)(g-i+1)}{2} +i -g$. Since our choices of the $Y_i$'s (and
$Y$) and hence our choices of the $Z_i$'s were arbitrary, we have
proved the second part of the Corollary as well.

To prove the third part, first observe that a dimension count (similar
to the case of $Y$) shows that the dimension of any irreducible
component $X$ of $E_q(V)$ is at least $\frac{q(q+1)}{2}
+\frac{(g-q)(g-q+1)}{2} -g$. Let $X'$ be the analogue of $Y'$ for $X$. Then
the images $X_q$ and $X_{g-q}$ of $X'$ by the two projections to
$\cA_q$ and $\cA_{g-q}$ are complete subvarieties of $\cA_q$ and
$\cA_{g-q}$ whose codimensions are at least $q$ and $g-q$ respectively
by part 1 of the Corollary. So we
have
\[ \begin{array}{c}
\frac{q(q+1)}{2} + \frac{(g-q)(g-q+1)}{2} -g \leq \hbox{dim}(X) =
\hbox{dim}(X') =
\hbox{dim}(X_q) + \hbox{dim}(X_{g-q}) \leq \\
\leq \frac{q(q+1)}{2} -q +
\frac{(g-q)(g-q+1)}{2} -(g-q) = \frac{q(q+1)}{2} +
\frac{(g-q)(g-q+1)}{2} -g \: .
\end{array} \]
Therefore we have equality everywhere and part 3 is proved.

Now let $V'$ be the analytic closure of $E_{1,g}(V)$ in $V$. Since, by
Theorem $1$, the set $E_1(V_{g-2})$ is dense in $V_{g-2}$ (which is a
curve), we see that $V'$ contains $Z_{g-2}$. Since all of our choices for
the $Y_i$ and $Y$ (and hence for the $Z_i$) were arbitrary, we see that
$V'$ contains $E_{1,g-2}(V)$. Repeating this reasoning, we see that $V'$
contains $E_{1,i}(V)$ for all $i$, hence $V'$ contains $E_1(V)$ and $V' =
V$ by Theorem $1$.  \hfill \qed

\vskip15pt
{\bf Proof of Corollary 2:}

Let $V$ be a complete codimension $g$ subvariety of $\tM_g'$ or
$\cA_g'$. Again, by Theorem $1$, the set $E_1(V)$ is nonempty. Let $Y
\subset V$ be an irreducible component of $E_1(V)$ and define $Y'$ as
in the proof of Corollary $1$. As in loc. cit. the variety $Y$ is a
complete subvariety of $V$, of codimension at most $g-1$ in $V$
(codimension exactly $g-1$ by Corollary $1$ if $V \subset \cA_g'$).

Suppose that $V \subset \tM_g'$. Again, since $Y'$ is irreducible and
complete and $\cA_1$ affine, the map $Y' \lra \cA_1$ is
constant, hence its differential has rank $0$ everywhere. It follows from
\cite{colpi} pages 172-173 that, for all $t \in Y \cap V_0$ and every
one-dimensional subvector space $W$ of $H^0(\omega_{C_t})$, the map $\mu :
W \otimes H^0(\omega_{C_t}) \lra T^*_t V$ is {\em not} injective. Since
this noninjectivity is a closed condition and $E_1(V_0)$ is dense in $V_0$,
it follows that it holds for all $t \in V_0$.

Therefore, by Proposition \ref{propNcurve} and with the notation there, for
all $t \in V_0$, there is a one-dimensional subvector
space $W_N$ of $H^0(\omega_{C_t})$ such that $N= \mu (W_N \otimes
H^0(\omega_{C_t}))$.

Let us globalize the constructions in the proof of Proposition
\ref{propNcurve}. Let $F_0$ be the Hodge bundle on $V_0$ and let $S^2 \bP
(F_0)$ be the quotient of the fiber product $\bP (F_0) \times_{V_0}
\bP (F_0)$ by the involution $\sigma$ exchanging the two factors of
the fiber product. Let $T^* \cM_g'$ be the cotangent bundle of
$\cM_g'$ and let $\cN_0 \subset T^* \cM_g' |_{V_0}$ be the conormal
bundle to $V_0$. Denote by $\cN''$ (resp. $\cN'$) the subcone of
decomposable tensors (resp. rank $2$ symmetric tensors) in $F_0
\otimes F_0$ (resp. $S^2 F_0$) lying in the inverse image of $\cN_0$
by the multiplication map $S^2F_0 \lra T^* \cM_g'
|_{V_0}$. Then, by Proposition
\ref{propNcurve} and with the notation there, the fibers of $\cN''$,
$\cN'$, and $\cN_0$ at $t$ are respectively $W_N \otimes
H^0(\omega_{C_t}) \cup H^0(\omega_{C_t}) \otimes W_N$, $\rho(W_N
\otimes H^0(\omega_{C_t}))$ and $\mu(W_N \otimes
H^0(\omega_{C_t}))$. Hence the morphism $m : \cN'
\lra \cN_0$ is an isomorphism because it is an isomorphism on each
fiber and the map $\bP(\cN'') \lra \bP(\cN_0)$ is a
double cover which splits on each fiber. Since the double cover of
$V_0$ parametrizing the rulings of the fibers of $\bP
(F_0) \times_{V_0} \bP (F_0)$ over $V_0$ is split, the double cover $\bP(\cN'')
\lra \bP(\cN') \cong \bP(\cN_0)$ is globally split and hence the
variety $\bP(\cN'')$ is the union of two subvarieties of $\bP (F_0)
\times_{V_0} \bP (F_0)$ exchanged by $\sigma$ and both isomorphic to
$\bP(\cN')$ (by the quotient morphism $\bP(F_0) \times_{V_0} \bP(F_0)
\lra S^2 \bP(F_0)$ ) and to $\bP(F_0)$ by either of the two projections
$\bP(F_0) \times_{V_0} \bP(F_0) \lra \bP(F_0)$. In particular, the two
components of $\bP(\cN'')$ are projective bundles on $V_0$ and $\cN''$
is the union of two vector bundles $\cN''_1$ and $\cN''_2$ with
respective fibers $W_N \otimes H^0(\omega_{C_t})$ and $H^0(\omega_{C_t})
\otimes W_N$ at $t$. Furthermore, we have $\cN''_1
\stackrel{\cong}{\lra} \cN_0 \stackrel{\cong}{\longleftarrow} \cN''_2$
(checked on fibers again). Since $\bP(\cN''_1)$ is isomorphic to $\bP
(F_0)$, there is a line bundle $\cW$ such that $\cN_1'' \cong \cW
\otimes F_0$. So $\cN_0 \cong \cW \otimes F_0$.

From the injection $\cN_1'' \inj F_0 \otimes F_0$ we deduce the
injection $\cW \inj F_0$ which is the composition of the morphism $\cW
\inj F_0 \otimes F_0 \otimes F_0^*$ (obtained from $\cW \otimes F_0 \cong
\cN''_1 \inj F_0 \otimes F_0$) with the morphism $F_0 \otimes (F_0 \otimes
F_0^*) \stackrel{id \otimes tr}{\lra} F_0$ which is the product of the
identity $F_0 \stackrel{id}{\lra} F_0$ and the trace morphism $F_0
\otimes F_0^* \cong End(F_0) \stackrel{tr}{\lra} \cO_{V_0}$.

For $V \subset \cA_g'$ the proof is similar to (and simpler than) the
above and uses Proposition \ref{propNppav} instead of Proposition
\ref{propNcurve}.
\hfill \qed

\section{Appendix: a remark on density in positive characteristic}
\label{sectpb}

In this section we use the notation of the introduction to denote moduli
spaces of curves and abelian varieties over an algebraically closed field
$k$ of characteristic $p > 0$. The subvariety $V_0$ of $\cA_g$
parametrizing ppav's of $p$-rank $0$ is a complete (connected if
$g > 1$ by \cite{oortnewton} (2.6)(c)) subvariety of codimension $g$ of
$\cA_g$ (see \cite{normoort}, (2) in the introduction and \cite{oortsubv},
the proof of Theorem 1.1a pages 98-99). We explain below how to deduce from
the results of \cite{demazure}, \cite{katz}, \cite{manin} and
\cite{oortnewton} that the moduli points of non-simple abelian varieties
are contained in a proper closed subset of $V_0$ when $g \geq 3$.

The formal group of an abelian variety $A$ of p-rank $0$ is isogenous to a sum
\[ \sum_{1 \leq i \leq r} G_{m_i,n_i} \]
where $m_i$ and $n_i$ are relatively prime positive integers for each $i$,
the sum $m_i+n_i$ is less than or equal to $g$ for all $i$, the formal
group $G_{m_i,n_i}$ has dimension $m_i$ and its dual is $G_{n_i,m_i}$ (see
\cite{manin} chapter IV, $\S 2$). The decomposition is symmetric,
i.e., the group $G_{m_i,n_i}$ appears as many times as
$G_{n_i,m_i}$. We call the unordered $r$-tuple $\left( (m_i,n_i)
\right)_{1 \leq i \leq r}$ the formal isogeny type
of the abelian variety. As in \cite{oortnewton}, we define the
Symmetric Newton Polygon of $A$ to be the lower convex polygon in the
plane $\bR^2$ which starts at $(0,0)$ and ends at $(2g,g)$, whose
break-points have integer coordinates and whose slopes (arranged in
increasing order because of lower convexity) are $\lambda_i =
\frac{n_i}{m_i+n_i}$ with multiplicity $m_i+n_i$ (i.e., on the
polygon, the $x$-coordinate grows by $m_i+n_i$ and the $y$-coordinate
grows by $n_i$). The polygon is symmetric in the sense that if the
slope $\lambda$ appears, then the slope $1- \lambda$ appears with the
same multiplicity. Following
\cite{oortnewton}, we shall say that the Newton Polygon $\beta$ is
above the Newton Polygon $\alpha$ if for all real numbers $x \in
[0,2g], y, z \in [0,g]$ such that $(x,z) \in
\beta$, $(x,y) \in \alpha$, we have $z \geq y$. We shall say that
$\beta$ is strictly above $\alpha$ if $\beta$ is above $\alpha$ and
$\beta \neq \alpha$. Again as in \cite{oortnewton}, for a Symmetric
Newton Polygon $\alpha$, we denote by $W_{\alpha}$ the set of points
in $\cA_g$ corresponding to abelian varieties whose Newton Polygon is
above $\alpha$. By \cite{demazure} page 91, Newton polygons go up
under specialization. By \cite{katz} page 143 Theorem 2.3.1 and
Corollary 2.3.2 (see also \cite{oortnewton}, 2.4), for any Newton
polygon $\alpha$, the set $W_{\alpha}$ is closed in $V_0$. By
\cite{oortnewton} Theorem (2.6)(a) and Remark (3.3), the abelian
variety $A_0$ with moduli point the generic point of $V_0$ has formal
isogeny type $((1,g-1),(g-1,1))$. Therefore, since $g \geq 3$, the
abelian variety $A_0$ is simple. Let $\alpha_0$ denote the Symmetric
Newton Polygon of $A_0$. The moduli point of a non-simple ppav of
$p$-rank $0$ is in $W_{\beta}$ for some Symmetric Newton Polygon
$\beta$ strictly above $\alpha_0$. Therefore the set of non-simple
ppav's in $V_0$ is contained in $\cup_{\beta \; strictly \; above \;
\alpha_0} W_{\beta}$. Since there are only a finite number of
Symmetric Newton Polygons (below the line $x=2y$ and) above
$\alpha_0$, we deduce that all points of $V_0$ corresponding to
nonsimple abelian varieties are in a proper closed subset of $V_0$
(which is $\cup_{\beta \; strictly \; above \; \alpha_0} W_{\beta}$).

Therefore $V_0$ is an example of a subvariety $V$
of codimension $g$ of $\cA_g$ (for all $g \geq 3$) or of $\tM_3$
such that $E_q(V)$ is not Zariski-dense in $V$ for any $q$.

\bibliographystyle{amsplain}

\providecommand{\bysame}{\leavevmode\hbox to3em{\hrulefill}\thinspace}

\end{document}